\documentclass[11pt]{article}

\usepackage{amssymb}

\usepackage{amscd}

\usepackage{amscd}

\title{Equivalences between cluster categories\thanks{Supported the NSF of China (Grants 10471071) and in part
 by Doctoral
Program Foundation of Institute of Higher Education(2003)}}

\author{Bin Zhu\thanks{E-mail: bzhu@math.tsinghua.edu.cn}}

\date{ \small {Department of Mathematical Sciences, \\ Tsinghua University,
100084
 Beijing, P. R. China \\  \vspace{0.4cm} Dedicated to Professor Claus M.
Ringel\\ on  the occasion of his sixtieth birthday.}}

\begin{document}

\maketitle

\def\s{\stackrel}

\def\gama{\gamma}
\def\Longrightarrow{{\longrightarrow}}
\def\P{{\cal P}}
\def\A{{\cal A}}
\def\F{\mathcal{F}}
\def\X{\mathcal{X}}
\def\T{\mathcal{T}}
\def\m{\textbf{ M}}
\def\t{{\tau }}
\def\b{\textbf{d}}
\def\K{{\cal K}}

\def\G{{\Gamma}}
\def\e{\mbox{exp}}

\def\righta{\rightarrow}

\def\s{\stackrel}

\def\ncong{\not\cong}

\def\mathbb{\NN}

\def\Hom{\mbox{Hom}}
\def\Ext{\mbox{Ext}}
\def\ind{\mbox{ind}}
\def\coprod{\amalg }
\def\L{\Lambda}

\newcommand{\uHom}{\operatorname{\underline{Hom}}\nolimits}
\newcommand{\End}{\operatorname{End}\nolimits}
\renewcommand{\r}{\operatorname{\underline{r}}\nolimits}
\def \text{\mbox}

%Egne definisjoner%

\begin{center}

\begin{minipage}{12cm}{\footnotesize\textbf{Abstract.}

\medskip

 Tilting theory in
cluster categories of hereditary algebras has been developed in
[BMRRT] and [BMR]. Some of them are already proved for hereditary
abelian categories there. In the present paper, all basic results
about tilting theory are generalized to cluster categories of
hereditary abelian categories. Furthermore, for any tilting object
$T$ in a hereditary abelian category $\mathcal{H}$, we verify that
the tilting functor Hom$_\mathcal{H}(T,-)$ induces a triangle
equivalence from the cluster category $\mathcal{C(H)}$ to the
cluster category $\mathcal{C}(A)$, where $A$ is the quasi-tilted
algebra End$_{\mathcal{H}}T.$ Under the condition that one of
derived categories of hereditary abelian categories $\mathcal{H},$
$\mathcal{H}'$ is triangle equivalent to the derived category of a
hereditary algebra, we prove that the cluster categories
$\mathcal{C(H)}$ and $\mathcal{C(H')}$ are triangle equivalent
 to each other if and only if $\mathcal{H}$ and $\mathcal{H}'$ are
derived equivalent, by using the precise relation between
cluster-tilted algebras (by definition, the endomorphism algebras of
tilting objects in cluster categories) and the corresponding
quasi-tilted algebras proved previously. As an application, we give
a realization of "truncated simple reflections" defined by
Fomin-Zelevinsky on the set of almost positive roots of the
corresponding type [FZ3, FZ4], by taking $\mathcal{H}$ to be the
representation category of a valued Dynkin quiver and $T$ a
BGP-tilting object(or APR-tilting, in other words).

\medskip

\textbf{Key words.} Tilting objects, cluster categories,
cluster-tilted algebras, BGP-reflection functors.

\textbf{Mathematics Subject Classification.} 16G20, 16G70, 19S99,
17B20.}
\end{minipage}
\end{center}

\newpage

\begin{center}

\textbf{1. Introduction}\end{center}

 Given a hereditary abelian category $\mathcal{H}$ with tilting objects,
 the orbit category
  of the (bounded)
  derived
 category $D^b(\mathcal{H})$ of $\mathcal{H}$ by
 its automorphism $F=[1]\tau^{-1}$ is again a triangulated category [Ke2],
 called cluster categories of type
 $\mathcal{H}$ and denoted simply by $\mathcal{C(H)}$. If $\mathcal{H}$
 is the category of representations of a Dynkin quiver $Q$, the
 corresponding cluster category $\mathcal{C}(Q)$ has been proved to be
 useful [BMRRT]: it provides a natural realization of
 clusters of corresponding cluster algebras, more precisely, there is
  an one to one correspondence
 from isoclasses of indecomposable objects in $\mathcal{C}(Q)$ to
 cluster variables
 of the corresponding cluster algebras. Under this correspondence,
 the basic tilting objects in $\mathcal{C}(Q)$ correspond
 to clusters (These results are generalized to all non-simply laced Dynkin types in [Z1], and recently are generalized to all acyclic quivers in [CK]).
 Clusters and cluster algebras are defined and studied by
 Fomin and Zelevinsky [FZ1-4][BFZ]. These algebras
 are
 defined so that it designs an algebraic framework for total positivity
  and canonical bases in semisimple algebraic groups.
    There are interesting connections to their theory in many directions
    [FZ1-4][BFZ] [CFZ] [GSV], amongst them to representation theory
    of quivers, in particular, to tilting theory [MRZ] [BMRRT] [BMR]
    [CFZ] [Z1].  Tilting theory in
    cluster category is really
 an  extension of classical titling theory of module category. The
  tilting objects in
 $\mathcal{C}(Q)$ are, on the one hand, corresponding to clusters  of corresponding cluster
 algebras (in
 simply-laced Dynkin type [BMRRT], in all Dynkin cases [Z1] and in all simply-laced cases [CK]); on the other hand, a generalization
 of tilting modules over hereditary algebras, for example, the
 endomorphism algebra of a tilting object in
 $\mathcal{C}(Q)$ may be self-injective.
\medskip

The aims of the paper are two-fold: The first one is to generalize
Buan-Marsh-Reiten theorem in [BMR] to the setting of hereditary
abelian categories with tilting objects. Buan-Marsh-Reiten theorem
says that the tilting functor Hom$_{\mathcal{C}(H)}(T, -)$ gives an
equivalence from the quotient
 $\mathcal{C}(H)/\mbox{add}\tau T$ of cluster category  to
 the module category of
 the cluster-tilted algebra of $T$. We prove that the same is true when
mod$H$ is replaced by any hereditary abelian category. Our proof for
the general result is obtained by a triangulated realization of that
 for Buan-Marsh-Reiten theorem, and simplifies the original proof.

  The second aim is to study the
triangle equivalences between cluster categories. We prove a "Morita
type" theorem for cluster categories. We verify the fact that any
 standard equivalence between two derived categories of hereditary
abelian categories induces a triangle equivalence between the two
corresponding cluster categories and prove that the inverse also
holds provided one of the hereditary abelian categories is derived
equivalent to a hereditary algebra. This first part is used in the
rest of paper and it was used in literatures, for examples: [BMRRT]
[BMR], and it is proved also in the updated version of [Ke2].  The reason why we verify the fact is that,
 the explicit expression of the triangle functor between triangulated orbit categories is
very interesting and useful, in particular, in the special case of
the triangle equivalences when it is induced by a
Bernstein-Gelfand-Ponomarev reflection tilting module or APR
tilting. It provides a realization of the "truncated simple
reflections" on the set of almost positive roots [Z1] (note that it
is proved in [Z2] that these triangle functors induce isomorphisms
of cluster algebras which are useful). By using this realization,
one can simplify some essential part of the quiver-theoretic
interpretation for generalized associahedra in the sense
 of Fomin-Zelevinsky [FZ4][CFZ] [MRZ] [Z1].
  In our proof of  "Morita type" Theorem for cluster categories, we use the
    precise relation between cluster-tilted algebras and the corresponding quasi-tilted
    algebras.
       Starting from a
basic tilting object $T$
 in $\mathcal{H}$ (note that "basic" means "multiplicity
free", we will assume that tilting objects are basic in the rest of
the paper), $T$ is a tilting object in cluster category
$\mathcal{C(H)}$ (note that any tilting object in $\mathcal{C(H)}$
 can be obtained from a tilting object in a hereditary abelian
 category $\mathcal{H}'$,
 derived equivalent to $\mathcal{H}$, so we don't loss the generality if we
 start from these tilting objects),
  then
the cluster-tilted algebra $\mbox{End}_{\mathcal{C(H)}}T$ of $T$ is
the trivial extension of the quasi-tilted algebra
$A=\mbox{End}_{\mathcal{H}}T$ of $T$ with the $A-$bimodule
Hom$_{\mathcal{H}}(T,\tau^2T)$. Some consequences follow: the
relations on Gabriel quivers and Auslander-Reiten quivers between
these two algebras become clear [BMR]. When $\mathcal{H}$ is the
module category of a hereditary algebra $H=KQ$, where $K$ is a
field, an important feather on
 the cluster category $\mathcal{C}(H)$ is that the generalization
 of APR-tilting at any
vertex is allowed, i.e., if $T'(i)$ is a projective module with all
indecomposable projective modules but one, say, $P(i)$, as its
direct summands . Then $T(i)=T'(i)\oplus \tau ^{-1}E(i)$ is a
tilting object in $\mathcal{C}(H)$. We assume $T(i)\in \mbox{mod}H$,
since $\tau^{-1}E(i) \notin \mbox{mod}H$ if and only if $E(i)$ is
injective, and in this case, we can replace $H$ by another
hereditary algebra $H'$, derived equivalent to $H$, with
$\tau^{-1}E'(i) \in \mbox{mod}H'.$ It follows from our result that
the cluster-tilted algebra $\mbox{End}_{\mathcal{C}(H)}T$ isomorphic
to End$_H(T)\ltimes D\mbox{Hom}_H(T, \tau E(i))$. When $i$ is a sink
or
 a source in $Q$, the cluster-tilted algebra goes back to the tilted
algebra of the same tilting module.

\medskip

This paper is organized as follows: in Section 2, some notions which
will be needed later on are recalled. Some basic properties of orbit
categories and cluster categories are given.  In particular, the
result that any standard equivalence between two derived categories
of hereditary abelian categories induces a triangle equivalence
between the corresponding orbit categories is verified; it is proved
that any almost complete tilting object in $\mathcal{C(H)}$ can be
completed to a tilting object in exactly two ways. In Section 3, we
prove that for a tilting object $T$ in a hereditary abelian category
$\mathcal{H}$, the cluster-tilted algebra
$\mbox{End}_{\mathcal{C(H)}}T$ is a trivial extension of a
quasi-tilted algebra with a bimodule. We explain through examples
how and what
 the Gabriel quivers or Auslander-Reiten quivers
  of the two algebras are related. We also prove the generalization
  of Buan-Marsh-Reiten
  theorem that the tilting
  functor Hom$_{\mathcal{C(H)}}(T,-)$ induces an equivalence from the
  quotient category $\mathcal{C}(H)/\mbox{add}\tau T$ of $\mathcal{C(H)}$ to the module category of
  cluster-tilted algebra. Our proof simplifies the original one in [BMR].
 In the final section, under the condition that one of the hereditary
 abelian categories is derived equivalent
 to a hereditary algebra, we prove that two cluster categories of
 hereditary abelian categories are triangle equivalent
 each other if and only if the two derived categories of hereditary
 abelian categories are triangle equivalent each other.
  As applications, we give a quiver
realization of "truncated simple reflections" on the set of almost
positive roots in all Dynkin types (simply laced  or non-simply
laced) and also give a quiver realization of Weyl generators of
Weyl group of any Kac-Moody Lie algebra.

 \begin{center}

\textbf{2. Basics on orbit categories and cluster categories.}
\end{center}

Let $\cal{H}$ be a hereditary abelian category with tilting objects and with finite dimensional Hom-spaces and Ext-spaces over a field
$K$, and denote by $\mathcal{D} =
    D^{b}(\cal{H})$ the bounded derived category of $\cal{H}$ with shift
functor $[1]$. For any category $\cal{E}$, we will denote by
$\ind\cal{E}$ the subcategory of isomorphism classes of
indecomposable objects in $\cal{E}$; depending on the context we
shall also use the same notation to denote the set of isomorphism
classes of indecomposable objects in $\cal{E}$. We write
$\mathcal{D} =
    D^{b}(\cal{H})$. Throughout the rest of paper, $D$ denotes the usual duality
Hom$_K(-,K)$ [ARS] [Rin].
\medskip

Let $G \colon \cal{D} \to \cal{D'}$ be a standard equivalence, i.e.,
$G$ is isomorphic to the derived tensor product
  $$ X\otimes _A-: D^b(A)\rightarrow D^b(A')$$
  for some complex $X$ of $A'-A-$bimodules.

 Following [Ke2],  we also assume $G$ satisfies the following properties:

\begin{itemize}
\item[(g1)]{For each $U$ in $\ind H$, only a finite number
of objects $G^n U$, where $n \in \mathbf{Z}$, lie in $\ind H$.}
\item[(g2)]{There is some $N \in \mathbf{N}$ such that
$\{U[n] \mid U \in \ind H, n \in [-N,N] \}$ contains a system of
representatives of the orbits of $G$ on $\ind \cal{D}$.}
\end{itemize}

We denote by $\mathcal{D}/ G$ the corresponding factor category.
The objects are by definition the $G$-orbits of objects in
$\cal{D}$, and the morphisms are given by
$$\Hom_{\mathcal{D}/G}(\widetilde{X},\widetilde{Y}) =
\oplus_{i \in \mathbf{Z}}
 \Hom_{\mathcal{D}}(X,G^iY).$$
Here $X$ and $Y$ are objects in $\cal{D}$, and $\widetilde{X}$ and
$\widetilde{Y}$ are the corresponding objects in $\mathcal{D}/G$
(although we shall sometimes write such objects simply as $X$ and
$Y$). The composition is defined in the natural way: if
$\tilde{f}: \widetilde{X}\rightarrow \widetilde{Y}$ and
$\tilde{g}: \widetilde{Y}\rightarrow \widetilde{Z}$ with $f:
X\rightarrow G^nY$ and $g: Y\rightarrow G^mZ$, then
$\tilde{f}\circ \tilde{g}$ is defined to be $\widetilde{fG^ng}$,
the image of composition of maps $f$ and $G^ng$ in $\mathcal{D}$.
The factor category $\mathcal{D}/G$ is Krull-Schmidt [BMRRT] and
is a triangulated [Ke2]. The canonical functor $\pi \colon
\mathcal{D} \longrightarrow \mathcal{D}/G: X\mapsto \widetilde{X}$
is a covering functor of triangulated categories [XZ2]. It sends
triangles to triangles. We remark that not all the triangles in
$\mathcal{D}/G$ are obtained
 as images of triangles in $\mathcal{D}$ under $\pi$. The shift in
 $\mathcal{D}/G$ is induced by the shift in
$\cal{D}$, and is also denoted by $[1]$. In both cases we write as
usual $\Hom(U,V[1]) = \Ext^1(U,V)$. We then have
$$\Ext^1_{\mathcal{D}/G}(\widetilde{X},\widetilde{Y}) =
\oplus_{i \in \mathbf{Z}} \Ext^1_{\mathcal{D}}( X, G^iY),$$ where
$X,Y$ are objects in $\cal{D}$ and $\widetilde{X},\widetilde{Y}$ are
the corresponding objects in $\mathcal{D}/G$. We shall mainly be
concerned with two special choices of functor  $F = \tau^{-1}[1]$ or
$F = [2]$ where $\tau$ is the Auslander-Reiten translation in
$\cal{D}$. In the first case, the factor category $\mathcal{D} /
\tau^{-1}[1]$ is called the cluster category of type $\mathcal{H}$,
which is denoted by $\mathcal{C}(\mathcal{H})$ (compare [BMRRT]). If
$\mathcal{H}$ is the module category of a hereditary algebra $H$ or
equivalently the category of representations of a valued quiver $Q$,
we denote the corresponding cluster category by $\mathcal{C}(H)$ or
$\mathcal{C}(Q)$ respectively. In the second case the factor
category $\mathcal{D} / [2]$ is called the root category of
 type $\cal{H}$, and we denote it by $\mathcal{R}(\mathcal{H})$
 (compare [H1] [H3] [XZZ]). When  $\mathcal{H}$ is
the module category of a hereditary algebra $H$ or a valued quiver
$Q$, we denote the corresponding root category by $\mathcal{R}(H)$
or $\mathcal{R}(Q)$ respectively.

  \medskip
Throughout the paper, $\cal{H}$ is assumed to be a hereditary
abelian category with titling objects. In this case, the
Grothendieck group $K_o(\cal{H})$ is a free abelian group of finite
rank. We recall that an object $T$ in $\cal{H}$
 is called a tilting object if Ext$^1_{\cal{H}}(T,T)=0$ and any
object $X$ with Ext$^1_{\cal{H}}(T,X)=\mbox{Hom}_{\cal{H}}(T,X)=0$
must be zero [HRS]. If $T$ is a tilting object in $\cal{H}$, then
the endomorphism algebra $A=\mbox{End}_{\mathcal{H}}(T)$ is called
a quasi-tilted algebra [HRS]. There are associated torsion pairs
$(\T, \F)$ in $\cal{H}$, and $(\mathcal{X},\mathcal{Y})$ in
$\mathcal{A}=$ mod$A$, such that there are equivalences of
categories Hom$_{\mathcal{H}}(T,-): \mathcal{T}\rightarrow
\mathcal{Y}$ and Ext$^1_{\mathcal{H}}(T,-): \mathcal{F}\rightarrow
\mathcal{X}$. In addition there is an induced equivalence of
derived categories
$$\mbox{RHom}(T,-): D^b(\mathcal{H})\rightarrow D^b(\mathcal{A}),$$
which is simply
denoted by $R(T)$ (compare [KZ] [H1]).

We recall the notation of Grothendieck groups of triangulated
categories $\A$ from [H1].
  Let $K(\A)$ be the free abelian
 group generated by representatives of the isomorphism classes of
 objects in $\A $.  The Grothendieck group $ K_{0}(\A )$ of $\A$ is the
 factor
group of $ K$ modulo the subgroup generated by elements of the
forms: $[A] + [C] - [B]$ corresponding to triangles $
A\longrightarrow B\longrightarrow C\longrightarrow TA.$

  \medskip

 \textbf{Proposition 2.1.} Let $\cal{H}$ be a hereditary abelian
 category with tilting objects
  and $G$ a triangle equivalence satisfying $(g_1),\ (g_2)$. Then
  $ K_{0}(\mathcal{D}/G
  )\\ \cong \textbf{Z}^n/H$,
  where $H$ is a subgroup of $\textbf{Z}^n$.
 \medskip

 \textbf{Proof.} We have a covering functor  $\pi :\
\mathcal{D}\rightarrow  \mathcal{D}/G $, which induces a surjective
group morphism
 $\pi _1: K(\mathcal{D})\longrightarrow  K(\mathcal{D} /G): [X]\mapsto
[\widetilde{X}],$ since $\pi$ is a triangle functor, $\pi _1$
induces a surjection from the Grothendieck group $K_0(\cal{D})$
 to the Grothendieck group $K_0(\mathcal{D}/G)$. It follows that
 $K_{0}(\mathcal{D}/G )\simeq K_{0}(\mathcal{D})/H$ and is isomorphic
 to $\textbf{Z}^n/H$ since from Theorem 4.6 in
Chapter I in [HRS] $K_{0}(\mathcal{D})\cong \textbf{Z}^n$ where $n$
is a positive integer. The proof is finished.
\medskip

\textbf{Remark 2.2.} We note that for some non-trivial orbit
triangulated categories, its Grothendieck groups might equal to
zero. For example $K_0(\mathcal{C}(Q))=0$ when $Q$ is a Dynkin
quivers of type $A_{2n}$. In contrary to the root category, the
Grothendieck group $K_0(\mathcal{R}(Q))$ is $\textbf{Z}^n$ for any
Dynkin type, where $n$ is the number of vertices of the quiver $Q$.

\medskip

 We recall the notation of exceptional set and of tilting set in
$\mathcal{C}(\cal{H}),$ from [BMRRT].
  A subset $B$ of ind$\mathcal{C}(\cal{H})$ is called exceptional
  if Ext$^1_{\mathcal{C}(\cal{H})}(X,Y)=0$ for any $X,\ Y\in
  B;$ An exceptional set $B$ is a tilting set if
  it is maximal with respect to this property. An object $T$ in
  $\mathcal{C}(\cal{H})$ is called
  tilting object if Ext$_{\mathcal{C}(\cal{H})}^1(T,T)=0$ and $T$
  has a maximal number of non-isomorphic
  direct summands. An object $M$ is called an almost complete tilting
  object if it is not a tilting object and there is an indecomposable
  object $X$ such that $M\oplus X$ is a tilting object. A subset $B$ of
  $\mathcal{C}(\cal{H})$ is a tilting set if and only if the direct
  sum of all objects in $B$ is a basic tilting object. For a tilting
  object $T$ in the cluster category
  of $\mathcal{H}$, The
  endomorphism ring End$_{\mathcal{C(H)}}T$ is called the
  cluster-tilted algebra of $T$.

  \medskip

Let $G \colon \cal{D} \to \cal{D}$  be a standard equivalence,
which is assumed to satisfy the properties (g1), (g2). Let $\Phi:
\mathcal{D}\longrightarrow \mathcal{D'}$ be a standard triangle
equivalence. We set $G'= \Phi G\Phi^{-1}$. Then $G'$ is a
 standard equivalence of $\mathcal{D}'$, it also satisfies the properties
(g1) and (g2). For most applications, we set $G=\tau ^m [n]$ for
some $m, \ n \in \textbf{Z},$ and then $G'$
 also equals to $\tau ^m [n]$.

\medskip

 \textbf{Definition 2.1.} We define the functor $\Phi_G$ from
 $\mathcal{D}/G$ to
$\mathcal{D}'/G'$ as follows: for $\widetilde{X}\in \mathcal{D}/G$
with $X\in \mathcal{D},$
 we set
$\Phi _G(\widetilde{X})=\widetilde{\Phi (X)}.$ For morphism
$\tilde{f}:\ \widetilde{X}\rightarrow \widetilde{Y}$, we set $\Phi
_G(\tilde{f}):\Phi_G(\widetilde{X}) \rightarrow \Phi_G(\widetilde{
Y})$ to be the map $\widetilde{\Phi (f)}.$
\medskip

The following result is proved in the updated version of
[Ke2], by using dg set-up. We will give some applications of it.

 \medskip

 \textbf{Proposition 2.3. } {\it Let $\Phi$ and $G$  be as above.
  Then $\Phi_G$
  is a triangle equivalence from $\mathcal{D}/G$ to
  $\mathcal{D}'/G'.$}

\medskip

\textbf{Proof.} First of all, we verify the definition is
well-defined: let $\widetilde{X}= \widetilde{Y}\in \mathcal{D}/G$
with $X, Y \in \mathcal{D}.$ Then we have that $Y=G^i(X)$ for some
integer $i$. It follows that $\Phi(Y)=\Phi G^i(X)\cong G'^i\Phi
(X)$. Hence $\Phi _G(\widetilde{X})= \Phi_G(\widetilde{Y})$. The
action of $\Phi_G$ on morphisms is induced by $\Phi$ on morphisms in
$\mathcal{D}$, and we have the commutative diagram as follows:

\[ \begin{CD}
 \oplus_{i\in \mathbf{Z}}\Hom _{\mathcal{D}}(X,G^iY)
  @>\Phi>>  \oplus_{i\in \mathbf{Z}}\Hom _{\mathcal{D}'}
  (\Phi ( X),
   G'^i\Phi(Y))\\
@V\wr VV  @VV\wr V  \\
\Hom _{\mathcal{D}/G}(\widetilde{X},\widetilde{Y}) @>\Phi_G
>> \Hom _{\mathcal{D}'/G'}(\Phi(\widetilde{X}),\Phi(\widetilde{Y}))
\end{CD} \]

 Therefore $\Phi_G$ is faithful and full. It is easy to see it is dense
 since $\Phi $ is equivalent. Then $\Phi_G$ is an equivalence from
 $\mathcal{D}/G$ to $\mathcal{D}'/G'$. Combining with that $\Phi _G$ is a triangle functor in section 9.4 in [Ke2], we have that
 $\Phi _G$ is a triangle equivalence from $\mathcal{D}/G$ to
  $\mathcal{D}'/G'.$ The proof is finished.

 \medskip

Under this triangle equivalence, tilting objects correspond to
tilting objects, and they have isomorphic endomorphism rings.
\medskip

 \textbf{Corollary 2.4.} {\it Let $\mathcal{H}$ (or $\mathcal{H}')$
 be a hereditary abelian category,
  $\Phi : \newline D^b(\mathcal{H}) \rightarrow D^b(\mathcal{H}')$
  a standard triangle equivalence and $G$ (or $G')$
  be as in Proposition 2.3. Let $T$ be an object in $\mathcal{C(H)}.$
   Then  $T$ is a tilting
   object in $\mathcal{C(H)} $ if and only so is $\Phi _G(T)$ in
   $\mathcal{C(H')}
   .$ Moreover End$_{\mathcal{C(H')}}(\Phi _G(T))\cong
   \mbox{End}_{\mathcal{C(H)}}(T)
   .$}

 \medskip

  If $\Phi$ is induced by a tilting object $T$ in $\cal{H}$, i.e.,
  $$\Phi =\mbox{RHom}(T,-): D^b(\mathcal{H})\rightarrow D^b(A),$$
  which is simply
denoted by $R(T)$, where $A$ is the endomorphism algebra of $T$,
then we have the following consequence:

\medskip

 \textbf{Corollary 2.5. }{\it Let $T$ be a tilting object in $\cal{H}$.
  Then $R_G(T)$
  is a triangle equivalence from $D^{b}(\mathcal{H})/G$ to
  $D^{b}(\mathcal{A})/G'.$}

\medskip

We will prove the converse of Proposition 2.3 when the orbit categories
are cluster categories and give some applications of Corollary 2.4.
in Section 4.
\medskip

In the rest of this section, we will prove some basic properties on
 tilting objects in a cluster category $C(\cal{H})$, where $\cal{H}$
is assumed a hereditary abelian category with tilting objects and
with Grothendieck group $\textbf{Z}^n$. These properties were
 proved in [BMRRT] when $\cal{H}$ is a module category of a
finite-dimensional algebra, and hold in the general case, which we
will show in the following.
\medskip

 \textbf{Proposition 2.6. }{(a). \it Let $T$ be a basic tilting
 object in $\mathcal{C}(\cal{H})$, where
 $\cal{H}$ is a hereditary abelian category with Grothendieck group
 $\textbf{Z}^n$. Then

  (i). $T$ is induced by a basic tilting object in a hereditary
  abelian category $\cal{H}'$, derived equivalent to
  $\cal{H}$.

  (ii). $T$ has $n$ indecomposable direct summands.

  (b). Any basic tilting objects in $\cal{H}$ induces a basic
  tilting objects for $\mathcal{C}(\cal{H})$.}

\medskip

\textbf{Proof.} The proof for statement (b) follows from the
definitions of tilting objects in various categories, namely, for a
tilting object $T$ in $\mathcal{H}$, we have that
Ext$^1_{\mathcal{C(H)}}(T,T)\cong
\mbox{Ext}^1_{\mathcal{H}}(T,T)\oplus
D\mbox{Ext}^1_{\mathcal{H}}(T,T)=0 $ and $T$ has $n$ non-isomorphic
indecomposable summands in $\mathcal{C(H)}$. Then  $T$ is a tilting
object in $\mathcal{C(H)}$. For the proof of (a), we note that any
hereditary abelian category with tilting objects is derived
equivalent to a module category of a hereditary algebra or to a
category of coherent sheaves over a weighted projective space [H2].
By Proposition 2.3. and Corollary 2.4, we can shift the proof to the
hereditary algebra case and the coherent sheaves case. For the first
case, all the statements were proved in [BMRRT]. For the second
case(compare [BMRRT]), through a suitable derived equivalence, we
may assume the hereditary abelian category $\cal{H}$ has no
projective or injective objects. Then the tilting objects in
$\mathcal{C}(\cal{H})$ and in $\cal{H}$ are 1-1 corresponding.
Therefore the statements in (a) hold.  The proof is finished.

\medskip

 \textbf{Proposition 2.7. }{\it Any exceptional object $\overline{T}$
 in $\mathcal{C}(\cal{H})$, where
 $\cal{H}$ is a hereditary abelian category with Grothendieck group
 $\textbf{Z}^n$, can be extended to a tilting object.
  If $\overline{T}$ is an almost complete basic tilting object, then
  $\overline{T}$ can be completed to a basic
tilting object in $\mathcal{C}(\cal{H})$ in exactly two different
ways.}

\medskip

\textbf{Proof.} We use the same strategy as Proposition 2.6. to
prove it. Since any hereditary abelian category with tilting objects
is derived equivalent to a module category of hereditary algebra or
to a category of coherent sheaves over a weighted projective space.
By the Proposition 2.3. and Corollary 2.4, we shift the statement to the
hereditary algebra case and the coherent sheaves case. For the first
case, all the statements are proved in [BMRRT]. For the second case,
through derived equivalence, we may assume the hereditary abelian
category $\cal{H}$ has no projective or injective objects. Then any
almost complete tilting object in $\mathcal{C}(\cal{H})$ are induced
from $\cal{H}$. By a result of Happel and Unger (Section 3 in [HU]),
there are exactly two complements of $\overline{T}$ in $\cal{H}$.
The proof is finished.

 \begin{center}

\textbf{3. Cluster-tilted algebras.}
\end{center}

  Since any tilting object $T$ in cluster category $\mathcal{C}(\cal{H})$
  can be
  obtained from a tilting object in a hereditary abelian category
  $\cal{H}'$, derived
  equivalent to $\cal{H}$, we may assume that,
  without loss the generality (compare Proposition 2.3. and Corollary 2.4.),
   $T$
  is a tilting object in $\cal{H}$, and then it is a tilting object in
  $\mathcal{C}(\cal{H})$.
  We have the quasi-tilted algebras $A=\mbox{End}_{\cal{H}}T$ and the
  cluster-tilted algebra $\Lambda=\mbox{End}_{\mathcal{C(H)}}T$. We will
  use $\mathcal{H}[k]$ to denote the full
  subcategory of $D^b(\mathcal{H})$ consisting of objects $X[k]$
  with $X\in \mathcal{H}.$
\medskip

It is easy to see that the quasi-tilted algebras are factor algebras
of cluster-tilted algebras from the definition of cluster
categories. In the following result, we explain cluster-tilted
algebras as the trivial-extensions of quasi-tilted algebras. This
explain is helpful for us to understand the relations on
combinatorics between these two algebras (compare [BMR]). We also
need this result to prove the Morita type Theorem for cluster
categories. We remind that  $D$ denotes the usual duality
Hom$_K(-,K)$ in the following.
\medskip

\textbf{Proposition 3.1.} Let $\cal{H}$ and $\mathcal{C}(\cal{H})$
be as above. Denoted by $A=\mbox{End}_{\cal{H}}T$ and
 $\Lambda=\mbox{End}_{\mathcal{C(H)}}T$. Then $\Lambda
 =A\ltimes D\mbox{Hom}_{\cal{H}}(T,\tau^2 T)$, i.e., $\Lambda$ is a
 trivial-extension of $A$ by the $A-$bimodule $D\mbox{Hom}_{\cal{H}}
 (T,\tau^2
 T)$.
 \medskip

Before we give the proof, we recall the notation of trivial
extensions (compare [ARS]). Given an algebra $A$ and  an
$A-$bimodule $M$, we define the algebra $A\ltimes M$ as follows, the
elements are pair $(a,m)$ with $a\in A$ and $m\in M$, addition is
componentwise and multiplication is given by
$(a,m)(b,n)=(ab,an+mb)$. It is easy to see if $A$ and $M$ are finite
dimensional, then $A\ltimes M$ is a finite dimensional algebra with
$M^2=0$.
\medskip

\textbf{Proof of Proposition.} From the definition of cluster-tilted
algebras, we have that $\L =\mbox{End}_{\mathcal{C(H)}}T=
\mbox{Hom}_{D^b(\mathcal{H})}(T,T)\oplus
\mbox{Hom}_{D^b(\mathcal{H})}(T,\tau ^{-1}T[1])=A\oplus
\mbox{Hom}_{D^b(\mathcal{H})}(T,\tau ^{-1}T[1]).$
  Where $\mbox{Hom}_{D^b(\mathcal{H})}(T,\tau ^{-1}T[1])$ is a
 natural left $A-$right
 $\mbox{End}_{D^b(\mathcal{H})}(\tau^{-1}T[1])-$module. Since
 $\tau^{-1}[1]$ is an automorphism of derived category
 $D^b(\mathcal{H})$, $\mbox{End}_{D^b(\mathcal{H})}(\tau^{-1}T[1])\cong A,$
 and
$\mbox{Hom}_{D^b(\mathcal{H})}(T,\tau ^{-1}T[1])$ is an
   $A-$bimodule.  It follows from the composition rule of morphisms
   in orbit category $\mathcal{C(H)}$ (compare [BG]) that $\Lambda$ is
   a trivial extension
   of $A$ with the $A-$bimodule $\mbox{Hom}_{D^b(\mathcal{H})}(
   T,\tau ^{-1}T[1]).$
  The remaining thing is to show \newline $\mbox{Hom}_{D^b(\mathcal{H})}
  (T,\tau ^{-1}T[1])
  \cong D\mbox{Hom}_{\mathcal{H}}(T,\tau ^2T)$ as $A-$bimodules.
  The first, we can view
  $\mbox{Hom}_{\cal{H}}(T,\tau^2
 T)$ is a natural $A-$bimodule as follows. It is a natural left $A-$right
  End$_{\cal{H}}(\tau^2
 T)-$bimodule. We assume $T=T_1\oplus T_2$ with $\tau ^2 T_2=0$ and assume
 $T_2$ is maximal with
 respect to this property. Then in derived category $D^b(\cal{H})$,
 object $\tau ^2T_2$ lies
 in the part $\mathcal{H}[-1]$ of degree $-1$, and $\tau ^2T_1$ lies
 in the part $\mathcal{H}[0]$ of degree $0$, hence
 Hom$_{D^b(\mathcal{H})}(\tau ^2T_1, \tau^2T_2)=0.$
Therefore we have the isomorphisms as follows:
$$\begin{array}{lll}\mbox{End}_{\cal{H}}(T)&\cong &
\mbox{End}_{D^b(\mathcal{H})}(\tau ^2 T_1\oplus \tau ^2 T_2)\\
&\cong &\big(\begin{array}{cc}\mbox{End}_{D^b(\cal{H})}\tau ^2T_1,&
\mbox{Hom}_{D^b(\cal{H})}(\tau ^2T_1,\tau ^2 T_2)\\
\mbox{Hom}_{D^b(\cal{H})}(\tau ^ 2 T_2,\tau
^2T_1),&\mbox{End}_{D^b(\cal{H})}\tau ^2 T_2\end{array}\big)\\
&\cong& \big(\begin{array}{cc}\mbox{End}_{\cal{H}}T_1,
 &0\\
\mbox{Hom}_{\cal{H}}(T_2,T_1),&\mbox{End}_{\cal{H}}T_2\end{array}\big).
\end{array}$$   Under these isomorphisms, any right
$\mbox{End}_{\mathcal{H}}\tau ^2T $ (which is
$\mbox{End}_{\mathcal{H}}\tau ^2T_1)-$module is a right
$\mbox{End}_{\mathcal{H}}T-$module. In the following, we will prove
that, as an $A-$bimodule, $\mbox{Hom}_{D^b(\mathcal{H})}(T,\tau
^{-1}T[1])$ is isomorphic to $D\mbox{Hom}_{\mathcal{H}}(T,\tau
^{2}T)$.

$$\begin{array}{lll}\mbox{Hom}_{D^b(\mathcal{H})}(T,
\tau^{-1}T[1])&\cong &\mbox{Ext}_{D^b(\mathcal{H})}^1(T, \tau^{-1}T)
\\
&\cong &D\mbox{Hom}_{D^b(\mathcal{H})}(\tau ^{-1}T, \tau^{-1}T)
\\
&\cong &D\mbox{Hom}_{D^b(\mathcal{H})}(T, \tau^{2}T)
\\
&\cong &D\mbox{Hom}_{\mathcal{H}}(T, \tau^{2}T).\end{array}$$
 The proof is finished.

 \medskip

We give some examples to illustrate the Proposition.
\medskip

\textbf{Example 3.2.} Let $Q$ be the quiver
$$\begin{array}{ccccc}1& \rightarrow& 2& \rightarrow& 3\end{array}
$$
and let $H = KQ$ be the path algebra, where $K$ is a field.

Let $T$ be the tilting module $T = T_1 \oplus T_2 \oplus T_3 =E_1
 \oplus P_1 \oplus P_3$, let $A = \mbox{End}_{H}(T)$ be the
corresponding tilted algebra and $\L =
\mbox{End}_{\mathcal{C}(H)}(T)$ the cluster-tilted algebra. We
notice that the tilted algebra $A$ is given by the quiver $Q$ also
with $\underline{r}^2 = 0$. Since $\tau ^2 T\cong P_3$ and
$D\mbox{Hom}_{H}(T,\tau ^2T)\cong D\mbox{Hom}_H(P_3, P_3)$ is an
one-dimensional space over $K$, it should contribute one arrow in
the quiver of $\L$. It is easy to check that there is an arrow
from vertex $3$ to $1$ in the quiver of $\L$ since the right
$A-$action on Hom$_H(T, \tau ^2 T_1)$ is provided by End$_H(T_1)$
via the isomorphisms indicated in the proof of Theorem 3.1.
 Therefore the quiver of $\L$ is the following:

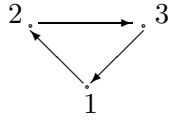
\begin{center}
 \setlength{\unitlength}{0.5cm}
 \begin{picture}(5,4)

\put(-0.6,3.2){2}\put(3.3,3.2){3} \put(1.4,0.8){1}
\put(1.4,1.6){\vector(-1,1){1.4}} \put(3,3){\vector(-1,-1){1.4}}

\put(0.2,3.2){\vector(1,0){2.5}}

\multiput(0,3.1)(3,0){2}{\circle{0.1}} \put(1.5,1.5){\circle{0.1}}

 \end{picture}
 \end{center}

with relations $\underline{r}^2=0$.

\medskip

\textbf{Example 3.3.} Let $Q$ be the quiver:
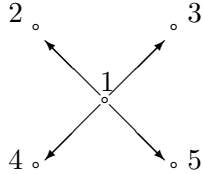
\begin{center}
 \setlength{\unitlength}{0.61cm}
 \begin{picture}(5,4)
 \put(-0.6,0){4}\put(3.3,0){5}
\put(-0.6,3.2){2}\put(3.3,3.2){3} \put(1.4,1,7){1}
\put(1.6,1.6){\vector(1,1){1.2}} \put(1.4,1.6){\vector(-1,1){1.2}}
\put(1.4,1.4){\vector(-1,-1){1.2}} \put(1.6,1.4){\vector(1,-1){1.2}}

\multiput(0,0.1)(3,0){2}{\circle{0.1}}
\multiput(0,3.1)(3,0){2}{\circle{0.1}} \put(1.5,1.5){\circle{0.1}}

 \end{picture}
 \end{center}
and $H=KQ$ the path algebra of $Q$.

Let $T$ be the tilting module $T = T_1 \oplus T_2 \oplus T_3\oplus
T_4 \oplus T_5 =R\oplus \tau^{-1}P_2\oplus \tau^{-1}P_3\oplus P_4
\oplus P_5$, where $R$ is the regular simple module with composition
factors $E_1, E_4, E_5$. Then the tilted algebra $A =
\mbox{End}_{H}(T)$ is the path algebra with relations: $ac-bf=0,
ec-df=0,$
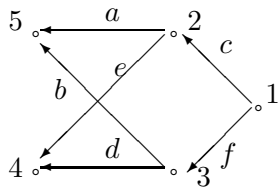
\begin{center}
 \setlength{\unitlength}{0.61cm}
 \begin{picture}(5,4)
 \put(-0.6,0){4}\put(3.5,-0.2){3}
\put(-0.6,3.2){5}\put(3.3,3.2){2} \put(5,1.6){1}
\put(4.7,1.6){\vector(-1,1){1.5}} \put(4.7,1.5){\vector(-1,-1){1.4}}
\multiput(2.8,0.2)(0,3){2}{\vector(-1,0){2.7}}

\put(2.8,0.2){\vector(-1,1){2.7}} \put(2.8,3.1){\vector(-1,-1){2.7}}

\multiput(0,0.1)(3,0){2}{\circle{0.1}}
\multiput(0,3.1)(3,0){2}{\circle{0.1}} \put(4.85,1.5){\circle{0.1}}

\put(1.5,3.4){$a$} \put(1.5,0.4){$d$}
\put(0.4,1.7){$b$}\put(1.7,2.2){$e$}\put(4,2.7){$c$}\put(4,0.3){$f$}

 \end{picture}
 \end{center}

  Since $\tau ^2 T\cong R$ and
$D\mbox{Hom}_{H}(T,\tau ^2T)\cong D\mbox{Hom}_H(T, R)$ is a
five-dimensional space over $K$. Any non-zero map from
$\tau^{-1}P_2$ (or $\tau^{-1}P_3$) factors through the identity map
of $R$, any non-zero map from $P_4$ (or $P_5$) factors through
$\tau^{-1}P_2$ and  through $\tau^{-1}P_3$. Then by Proposition
 3.1., to get the quiver of the corresponding cluster-tilted algebra
$\L$, we add two arrows in the quiver of $A$, one is from $4$ to
$1$, another from $5$ to $1$, with the additional relations:
$ga-he=0, hd-gb=0, cg=fg=ch=fh=0,$ i.e. $\L$ is the path algebra
with relations: $ac-bf=0, ec-df=0, ga-he=0, hd-gb=0, cg=fg=ch=fh=0,$

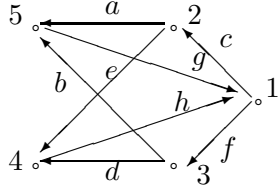
\begin{center}
 \setlength{\unitlength}{0.61cm}
 \begin{picture}(5,4)
 \put(-0.6,0){4}\put(3.5,-0.3){3}
\put(-0.6,3.2){5}\put(3.3,3.2){2} \put(5,1.6){1}
\put(4.7,1.6){\vector(-1,1){1.5}}
\put(4.7,1.5){\vector(-1,-1){1.4}}
\multiput(2.8,0.2)(0,3){2}{\vector(-1,0){2.7}}

\put(2.8,0.2){\vector(-1,1){2.7}} \put(2.8,3.1){\vector(-1,-1){2.7}}
\put(0.1,0.2){\vector(3,1){4.2}} \put(0.2,3.1){\vector(3,-1){4.2}}
\multiput(0,0.1)(3,0){2}{\circle{0.1}}
\multiput(0,3.1)(3,0){2}{\circle{0.1}} \put(4.85,1.5){\circle{0.1}}

\put(1.5,3.4){$a$} \put(1.5,-0.2){$d$}
\put(0.4,1.7){$b$}\put(1.5,2){$e$}\put(4,2.7){$c$}\put(4,0.3){$f$}

\put(3,1.3){$h$} \put(3.4,2.3){$g$}

 \end{picture}
 \end{center}

\textbf{Remark 3.4.} Since $\L$ is a trivial extension of $A$,
mod$A$ is embedded in mod$\L$ as the full subcategory consisting of
$\L-$modules $X$ which are annihilated by the idea
$D\mbox{Hom}_{\cal{H}}(T, \tau ^2 T)$. If we consider
Auslander-Reiten quivers of these algebras, we can get AR-quiver of
$A$ from that of $\L$ by deleting those vertices from which there
are non-zero maps (non-zero path) to $\tau ^2 T$ after identifying
the AR-quivers between $\L$ and $\mathcal{C(H)}/(\mbox{add}\tau T)$
[BMR]. For example, AR-quiver of $A$ in Example 3.2. is obtained
from AR-quiver of $\L$ by deleting $P_3$, since there is only one
indecomposable object, namely, $P_3$, from which there exist
non-zero map to $\tau ^2T\cong P_3.$

 \medskip

Let $H$ (or $H'(k)$) be the tensor algebra of an valued quiver
$(\m, \G, \Omega)$\newline  (resp. $(\m, \G, s_k\Omega)$, resp.)
with indecomposable projective modules $P_i, \ i\in \G _0.$ For
any $k$, let $T'(k)=\oplus _{i\not= k}P_j$. It is an almost
complete tilting object in $\mathcal{C}(H)$, and any almost
complete tilting object in $\mathcal{C}(H)$ can be obtained from
an almost complete tilting module over a hereditary algebra $H'$,
derived equivalent to $H$. There are exactly two ways to complete
$T'(k)$ into a tilting object: one is to plus $P_k$, another one
is to plus $\tau ^{-1}E_k$. The tilting object $T(k)=T'(k)\oplus
\tau ^{-1}E_k$ is called APR-tilting object in $\mathcal{C}(H)$ at
vertex $k$ in [BMR]. In case $k$ is sink in the quiver of $H$,
$T(k)$ is the usual BGP-tilting module, in this case $E_k=P_k.$
Denoted by $A(k)$ the tilted algebra of APR-tilting module $T(k)$.
The cluster-titled algebras, denoted by $\L (k)$, of APR-tilting
objects are described explicitly in the following way:
\medskip

 \textbf{Corollary 3.5.} Let $T(k)$ be the
APR-tilting object in
 $\mathcal{C}(H)$. Then if $k$ is sink or source in $\Omega$, $\Lambda
  (k)\cong H'(k)$; otherwise,
 $\L(k)=A(k)\ltimes D\mbox{Hom}_H(T'(k),\tau
 E_k)$.
\medskip

 \textbf{Proof.} We note that $k$ is sink if and only if $P_k$ is simple
 module $E_k$. In this case,
 $\tau ^2T(k)\cong 0$. It follows from Proposition 3.1 that $\Lambda (k)
 \cong H'(k).$ We also note that
 $k$ is source if and only if $E_k$ is injective, hence $\tau ^{-1}E_k$
 does not
 exist
  in $H-$module, but exists in $\mathcal{C}(H)$. If we consider the
  reflection at vertex $k$, we have that
  the hereditary algebra $H'(k)$ is derived equivalent to $H$.
  It follows that $R(S_k^-)(T(k))=H'(k).$ Then by Corollary 2.4., we
  have $\Lambda \cong H'(k).$
  Suppose $k$ is neither a sink nor a source.  Since Hom$_H(\tau ^{-1}E_k,
  \tau E_k)=0$, we have that
   $\mbox{Hom}_H(T(k), \tau ^2(\tau ^{-1}E_k))\cong
 \mbox{Hom}_H(T(k), \tau E_k)\cong \mbox{Hom}_H(T'(k), \tau E_k)$.
 The proof is finished.

 \medskip

 The cluster-tilted algebras of APR-tilting objects may be quite
 different from
 tilted algebras of APR-tilting modules. For example, the tilting
 object in Example
 3.2. is an APR-tilting at vertex $2$. The corresponding cluster
 tilted algebra is
 self-injective.
 \medskip

 From Corollary 3.5., one can easily determine all cluster-tilted
  algebras of APR-tilting objects in
 a given cluster category $\mathcal{C}(H)$.
\medskip

\textbf{Example 3.6.} Let $H=KQ$, where $Q$ is the quiver:
$$\begin{array}{ccccccc}&&&&5&&\\
&&&&\downarrow&&\\
1&\rightarrow&2&\rightarrow&3&\rightarrow&4.\end{array}$$ Now
vertices $1$ and $5$ are sources, then the cluster-tilted algebras
$\Lambda (1)$ and $\Lambda (5)$ corresponding to APR-tilting
objects are $H'(1)$ or $H'(5)$ respectively, where $H'(1)$ is the
quiver algebra $K(s_1Q)$ and $H'(5)$ is the quiver algebra
$K(s_5Q).$ Vertex $4$ is sink, then $\Lambda (4)$ is the quiver
algebra $K(s_4Q).$ We compute other cluster-tilted algebras
corresponding to APR-tilting objects at vertices $2$ or $3$ in
$\mathcal{C}(H).$  The APR-tilting module at vertex $2$ is
$T(2)=P_1\oplus \tau ^{-1}E_2\oplus P_3 \oplus P_4\oplus P_5 .$
The corresponding tilted algebra $A(2)$ is the following:
\begin{center}
 \setlength{\unitlength}{0.61cm}
 \begin{picture}(10,3)
 \put(1,0){2\vector(1,0){1.5}}
 \put(2.9,0){1\vector(1,0){1.5}}
 \put(4.8,0){3\vector(1,0){1.5}}
 \put(6.5,0){4, }
 \put(4.6,2.1){5\vector(0,-1){1.5}}
\put(2,0.2){$\alpha$} \put(3.8,0.2){$\beta$} \put(5.7,0.2){$\gamma$}
\put(4.3,1.3){$\rho$} \put(8.5, 0){$\beta \alpha =0$}

 \end{picture}
 \end{center}

Since $\tau E_2=\tau^{-1}P_1$, Hom$_H(T(2), \tau E_2)\cong
\mbox{Hom}_H(P_3\oplus P_5, \tau E_2)$ is two-dimensional, and the
non-zero morphism from $P_3$ to $\tau E_2$ factors through $P_5$.
Then to get the quiver of the corresponding cluster-tilted algebra
$\Lambda (2)$, we add
 an arrow $\delta$ from vertex $3$ to $2$ in the quiver of $A(2)$ with the
addition relations $\delta \beta=0$, $ \alpha \delta=0$.  That is

\begin{center}
 \setlength{\unitlength}{0.61cm}
 \begin{picture}(10,3)
 \put(1,0){2\vector(1,0){1.5}}
 \put(2.9,0){1\vector(1,0){1.5}}
 \put(4.8,0){3\vector(1,0){1.5}}
 \put(6.5,0){4, }
 \put(4.5,2.1){5\vector(0,-1){1.5}}
\put(2,0.2){$\alpha$} \put(3.8,0.2){$\beta$} \put(5.7,0.2){$\gamma$}
\put(4.3,1){$\rho$} \put(8.5, 0){$\beta \alpha=\delta \beta=\alpha
\delta=0$} \put(4.8,-0.2){\vector(-1,0){3.7}} \put(3,-0.6){$\delta$}
 \end{picture}
 \end{center}

The APR-tilting module at vertex $3$ is $T(3)=P_1\oplus P_2 \oplus
\tau ^{-1}E_3 \oplus P_4\oplus P_5.$ The corresponding tilted
algebra $A(2)$ is the following:
\begin{center}
 \setlength{\unitlength}{0.61cm}
 \begin{picture}(10,3)
 \put(1,0){2\vector(1,0){1.5}}
 \put(2.9,0){1\vector(1,0){1.5}}
 \put(3,2.5){3\vector(1,0){1.5}5}
 \put(5,2.5){\vector(0,-1){1.7}}
 \put(4.8,0){4, }
 \put(3,2.5){\vector(0,-1){1.7}}

\put(2,0.2){$\alpha$} \put(3.8,0.2){$\beta$} \put(4,2.7){$\gamma$}
\put(5.3,1.2){$\rho$}\put(2.7, 1.2){$\delta$} \put(8.5, 0){$\beta
\delta=\rho \gamma$}

 \end{picture}
 \end{center}

Since $\tau E_3=P_4$, Hom$_H(T(2), \tau E_3)\cong \mbox{Hom}_H(P_4,
\tau E_3)$ is one-dimensional. Then to get the quiver of the
corresponding cluster-tilted algebra $\Lambda (3)$, we add
 an arrow $\xi$ from vertex $4$ to $3$ in the quiver of $A(3)$ with the
addition relations $\xi \beta=\xi\rho=\delta \xi=\gamma\xi =0$. That
is

\begin{center}
 \setlength{\unitlength}{0.61cm}
 \begin{picture}(10,3)
 \put(1,0){2\vector(1,0){1.5}}
 \put(2.9,0){1\vector(1,0){1.5}}
 \put(3,2){3\vector(1,0){1.5}5}
 \put(5,2){\vector(0,-1){1.4}}
 \put(4.8,0){4, }
 \put(3,2){\vector(0,-1){1.4}}
\put(4.7,0.2){\vector(-1,1){1.5}} \put(4,1){$\xi$}
\put(2,-0.4){$\alpha$} \put(3.8,-0.4){$\beta$} \put(4,2.3){$\gamma$}
\put(5.3,1){$\rho$}\put(2.7, 1){$\delta$} \put(6, 0){$\beta
\delta=\rho \gamma$, $\xi \beta=\xi\rho=\delta \xi=\gamma\xi =0$}

 \end{picture}
 \end{center}

\medskip

  One of important results on cluster-tilted algebras is Theorem 2.2
  in [BMR] which gives precise relation between the cluster
  categories and module category over cluster-tilted algebras. We
  will generalize the result to the general setting, where the module
  categories over
  hereditary algebras are replaced by any hereditary abelian categories
  with tilting objects.
  Our proof simplifies the proof in [BMR] and uses approximations on
  triangulated
  categories.
  \medskip

  Let $T$ be a tilting object in $\cal{H}$. Then Hom functor
 $G=\mbox{Hom}_{\mathcal{H}}(T,-)$ induces a dense and full functor
 from the
  cluster category $\mathcal{C(H)}$ to $\L-$mod. Where $\L$ is the
  cluster-tilted algebra End$_{\mathcal{C(H)}}T$, the density and
  fullness of the functor is obtained from Proposition 2.1. in [BMR]
  (The proof there
  are not involved using hereditary algebra, so the proof also
  works in this general case). Since $G(\tau T)=\mbox{Hom}_{\mathcal{C(H)}}
  (T,\tau T)=0,$ there is an
  induced functor $\bar{G}: \mathcal{C(H)}/\mbox{add}(\tau
  T)\rightarrow \mbox{mod}\L$. We will prove that $\bar{G}$ is faithful
  in the following. Firstly we recall
  that there is an embedding which identifies
  $\mathcal{H}$ with the full subcategory of $D^b(\mathcal{H})$
  consisting of complexes which have zero components of any non-zero
  degree. We remaind the reader that $F$ denotes $\tau ^{-1}[1].$
\medskip

\medskip

\textbf{Theorem 3.7.} Let $T$ be a tilting object in $\cal{H}$ and
 $\Lambda=\mbox{End}_{\mathcal{C(H)}}T$ the cluster-tilted algebra.
 Then  $\overline{G}: \mathcal{C(H)}/\mbox{add}(\tau
  T)\rightarrow \mbox{mod}\L$ is an equivalence.

\medskip

\textbf{Proof.}  We only need to show that $\overline{G}$ is
faithful. Let $\bar{f}: \overline{M}\rightarrow \overline{N} $ be a
map between indecomposable objects in $\mathcal{C(H)}$. We can
assume that $\bar{f}$ is induced from a map $f: M\rightarrow N $ in
D$^b(\mathcal{H})$ with $M, \ N$ in $\mathcal{H}$ or
 $\mathcal{H}[1]$. We will show if
$\overline{G}(\bar{f})=0$, then $f$ factors through add$(\tau T)$.
It follows from $\overline{G}(\bar{f})=0$ that
Hom$_{D^b(\mathcal{H})}(T,f)=0$ and
Hom$_{D^b(\mathcal{H})}(T,Ff)=0.$

(I). We consider the minimal right add$T-$approximation of $M,$ $
M_2\s{\alpha_2}{\longrightarrow} M$ with $M_2\in \mbox{add}T$ which
exists since add$T$ contains finitely many indecomposable objects.
 Then we have
the triangle: $$(*)\ \ \ M_2\s{\alpha_2}{\longrightarrow}M\s{\alpha
_1}{\longrightarrow} M_1\longrightarrow M_2[1].$$
  We also have the triangle with $\beta _2: N{\rightarrow }N_2$ is
  the minimal left add$\tau ^2T-$approximation:
  $$(**)\ \ \ N_1\s{\beta _1}{\longrightarrow} N\s{\beta _2}{\longrightarrow}
N_2\longrightarrow N_1[1]$$ where $N_2 \in \mbox{add}\tau^2 T$. It
is easy to see $M_1\in \mathcal{H}\cup \mathcal{H}[1]$ and  $N_1\in
\mathcal{H}\cup \mathcal{H}[-1].$

We will prove that there exists a commutative diagram:

$$\begin{array}{rcl}M& \s{f}{\longrightarrow}&N\\
\alpha\downarrow&&\uparrow\beta   \\
M'& \s{g_{1}}{\longrightarrow}&N'
\end{array}$$

 with $M'_1, N'_1 \in \mathcal{H}$ and $M'_1$(or $ N'_1$) is direct
  summand of $M_1(N_1$, resp.)

 (1). It follows from Hom$_{D^b(\mathcal{H})}(T,Ff)=0$ that
 Hom$_{D^b(\mathcal{H})}(\tau^2T,f)=0$. It implies
 that $f \beta _2 =0$ since $N_2\in \mbox{add}\tau ^2T$. Then
  there is a map $f_1: M\longrightarrow N_1$
 with $f=f_1\beta _1.$

 (2). We prove that Hom$(T,f_1)=0$. By applying Hom$(T,-)$ to
 the triangle $(**)$, we have the exact sequence:
 $$ \mbox{Hom}(T,N_2[-1])\longrightarrow  \mbox{Hom}(T,N_1)
 \s{\mbox{Hom}(T,\beta _1)}{\longrightarrow}
  \mbox{Hom}(T,N)\longrightarrow  \mbox{Hom}(T,N_2).$$
 Since $ \mbox{Hom}(T,N_2[-1])\subseteq \mbox{Hom}(T,(\tau ^2T[-1])^m)
 =(D\mbox{Ext}^1(\tau T, T[1])^m=\newline(D\mbox{Hom}(\tau T, T[2]))^m =0,$
 for some positive integer $m$, the map $\mbox{Hom}(T,\beta _1)$ is mono. It follows that
 Hom$(T,f_1)$ is zero since the its composition with Hom$(T,\beta_1)$
 is Hom$(T,f)=0.$

 (3). By (2), we have $\alpha _2 f_1=0$. It follows that there
 exists a map $g: M_1\longrightarrow N_1$ with $f_1=\alpha _1 g$.

 (4). We write $M_1=M'_1\oplus M'_2$ and $N_1=N'_1\oplus N'_2$ with
 $M'_1$ and  $N'_1$ are maximal
 direct summand in $\mathcal{H}$ of $M_1$ and $N_1$ respectively. Then
 $g =(\begin{array}{cc}g_1, &o\\
 0&0\end{array}).$ Let $\alpha$ be the component of $\alpha _1$ on
 $M_1'$ and  $\beta$ the component of $\beta _1$
 on $N_1'.$ Then $f=\alpha g_1\beta$. Then we have proved there exists
 a commutative diagram
 which we proposed above. For simplicity, we assume that both $M_1$ and
 $N_1$ are in
 $\mathcal{H}$.

\vspace{0.3cm}

(II). We will prove that map $g_1: M_1\longrightarrow N_1$ factors
through
 add$\tau T$.
\medskip

By applying Hom$(T,-)$ to the triangle $(*)$, we can get exact
sequence \newline $\mbox{Hom}(T,M_2)\s{\mbox{Hom}(T,\alpha
_2)}{\Longrightarrow} \mbox{Hom}(T,M){\longrightarrow}
\mbox{Hom}(T,M_1)\longrightarrow \mbox{Hom}(T,M_2[1])=0,$ where
$\mbox{Hom}(T,\alpha _2)$ is surjective. It follows that
$\mbox{Hom}_{\mathcal{H}}(T, M_1)=0$, with the conditions that $T,\
M_1\in \mathcal{H}$. Then there is an embedding of $M_1$ into $\tau
T$. Let $0\longrightarrow
M_1\s{\gama}{\longrightarrow}E_1\longrightarrow E_2\longrightarrow
 0$ be an exact sequence in $\mathcal{H}$ with
 $\gama$ being the minimal left add$(\tau T)-$approximation of $M_1$,
 where $E_1\in \mbox{add}(\tau T).$
 Then $(***)\ \ \ \ M_1\s{\gama}{\longrightarrow}E_1\longrightarrow
 E_2\longrightarrow
M_1[1]$ is a triangle with
 $\gama $ is the minimal left add$(\tau T)-$approximation of $M_1$ in
 $D^b(\mathcal{H}).$ Of course we have that $E_1, \ E_2\in
 \mathcal{H}$

\medskip

 (1). To prove Hom$(\tau ^{-1}N_1,\tau T)\cong \mbox{Hom}(N_1,\tau ^2T)=0$.
 \medskip

 By applying Hom$(-,\tau ^2T)$ to the triangle $(**)$, we have the exact
 sequence with $\mbox{Hom}(\beta _2,\tau^2
 T)$ being surjective: $$\begin{array}{l}\mbox{Hom}(N_2,\tau^2T)
 \s{\mbox{Hom}(\beta _2,\tau^2 T)}{\longrightarrow}
 \mbox{Hom}(N,\tau^2T){\longrightarrow}\mbox{Hom}(N_1,\tau^2T)\\
 \longrightarrow
 \mbox{Hom}(N_2[-1],\tau^2T).\end{array}$$  But $\mbox{Hom}(N_2[-1],
 \tau ^2T)\subseteq \mbox{Hom}((\tau ^2T[-1])^m,\tau
 ^2T)=0.$ It follows that \\
 $\mbox{Hom}(N_1,\tau ^2T)=0$.

\medskip

(2). To prove that $\mbox{Hom}(T,E_2)=0$.
\medskip

 By applying Hom$(-,\tau T)$ to the triangle $(***)$, we have the
 exact sequence with surjective map $\mbox{Hom}(\gamma,\tau T)$:
   $$\mbox{Hom}(E_1,\tau T)\s{\mbox{Hom}(\gamma,\tau T)}{\longrightarrow}
    \mbox{Hom}(M_1,\tau T)\longrightarrow
  \mbox{Hom}(E_2[-1],\tau T)\longrightarrow\mbox{Hom}(E_1[-1],\tau T).$$
$\mbox{Hom}(E_1[-1],\tau T)\cong D\mbox{Ext}^1(T, E_1[-1])
\subseteq D\mbox{Ext}^1(T, (\tau T[-1])^n)\cong (D\mbox{Hom}(T,
\tau T))^n\\ \cong (\mbox{Ext}^1(T, T))^n=0.$  Therefore
Hom$(E_2[-1],\tau T)=0$ and then Hom$(T,E_2)=0.$ Therefore
$\mbox{Hom}_{\mathcal{H}}(T, E_2)=0.$
\medskip

(3). To prove $\mbox{Hom}(\tau ^{-1}N_1,E_2)=0$.

\medskip

 We
have that there is an embedding of $E_2$ into $\tau T$. Let
$0\longrightarrow E_2\s{\sigma}{\longrightarrow}X_1\longrightarrow
X_2\longrightarrow
 0$ be an exact sequence in $\mathcal{H}$ with
 $\sigma$ being the minimal left add$(\tau T)-$approximation of $E_2$,
 where $X_1\in \mbox{add}(\tau T).$
 Then $ E_2\s{\sigma}{\longrightarrow}X_1\longrightarrow
X_2\longrightarrow
 E_2[1]$ is a triangle with
 $\sigma$ being the minimal left add$(\tau T)-$approximation of $E_2$ in
 $D^b(\mathcal{H}),$ where $X_1, \ X_2\in
 \mathcal{H}$.
By applying Hom$(\tau ^{-1}N_1,-)$ to this triangle, we have the
exact sequence:$$\mbox{Hom}(\tau^{-1}N_1, X_2[-1])\longrightarrow
\mbox{Hom}(\tau^{-1}N_1,
E_2)\longrightarrow\mbox{Hom}(\tau^{-1}N_1, X_1).$$ Where
$\mbox{Hom}(\tau^{-1}N_1, X_1)\subseteq \mbox{Hom}(\tau^{-1}N_1,
(\tau T)^t)=0$ (by (1) in II).  We also have that
$\mbox{Hom}(\tau^{-1}N_1, X_2[-1])=\mbox{Hom}(\tau^{-1}N_1, \tau
X_2[-1])=D\mbox{Ext}^1(X_2[-1],N_1)=D\mbox{Ext}^2(X_2,N_1)=D\mbox{Hom}
(X_2,N_1[2])=0$
 (by $X_2, \ N_1 \in \mathcal{H}).$
  It follows that $\mbox{Hom}(\tau ^{-1}N_1,E_2)=0$ .
\medskip

  (4). By applying Hom$(-,N_1)$ to triangle $(***)$, we get the
  exact sequence:
  $$ \mbox{Hom}(E_2,N_1)\longrightarrow \mbox{Hom}(E_1,N_1)
  \longrightarrow\mbox{Hom}(M_1,N_1)\longrightarrow
  \mbox{Hom}(E_2[-1],N_1).$$
 Where $$\mbox{Hom}(E_2[-1],N_1)\cong\mbox{Ext}^1(E_2,N_1)\cong
 D\mbox{Hom}(N_1,\tau E_2)\cong
 \mbox{Hom}(\tau ^{-1}N_1, E_2)=0.$$ Therefore the map $g_1$ factors
 through $E_1\in \mbox{add}\tau T.$  The proof is finished.

 \medskip

\begin{center}

\textbf{4. Equivalences between cluster categories.}
\end{center}

\medskip

In Section 2,  we verified the fact that any triangle equivalence
between derived categories of hereditary categories induces a
triangle equivalence between corresponding triangulated orbit categories of derived categories by suitable automorphisms. Then
 standard equivalences induce triangle equivalences of cluster
categories. We will first proved the converse also holds in this
section (The problem whether the converse holds is suggested by
Professor Steffen Koenig, we thank him very much!), then
 we give some useful consequences on cluster categories
and root categories as applications.

\medskip

\textbf{Theorem 4.1.} Let  $\mathcal{H}_1$ and $\mathcal{H}_2$ be
hereditary abelian categories and one of them derived equivalent
to module category of a hereditary algebra. Then
$\mathcal{C}(\mathcal{H}_1)$ is triangle equivalent to
$\mathcal{C(H}_2)$ if and only if $\mathcal{H}_1$ is derived
equivalent to $\mathcal{H}_2.$
\medskip

\textbf{Proof.} The sufficiency is the special case of Proposition
2.3. in which $G=\tau ^{-1}[1]$ by using the "Morita Theorem" on derived categories
(compare [Ri] [Ke1]). We need to prove the necessity.
Let $\alpha : \mathcal{C}(\mathcal{H}_1)\longrightarrow
\mathcal{C}(\mathcal{H}_2)$ be a triangle equivalence. Suppose $T$
is a tilting object in $\mathcal{C}(\mathcal{H}_1).$ We will prove
that $\alpha (T)$ is also a tilting object in
$\mathcal{C}(\mathcal{H}_2).$  Since
Ext$_{\mathcal{C}(\mathcal{H}_2)}^1(\alpha (T), \alpha
(T))=\mbox{Ext}_{\mathcal{C}(\mathcal{H}_1)}^1(T,T)=0,$ $\alpha
(T)$ is an exceptional object in $\mathcal{C}(\mathcal{H}_2)$. It
follows from Theorem 2.7. or Proposition 3.2 in [BMRRT] that there
is an object $M$ such that $\alpha (T)\oplus M$ is a tilting
object in $\mathcal{C}(\mathcal{H}_2)$. Then $\alpha ^{-1}(\alpha
(T)\oplus M)\cong T \oplus \alpha ^{-1}( M)$ is an exceptional
object in $\mathcal{C}(\mathcal{H}_1)$. It follows from that $T$
is a tilting object that $\alpha ^{-1}( M)\in \mbox{add}T.$ Then
$M\in \mbox{add}\alpha(T).$ Hence $\alpha (T)$ is a tilting object
in $\mathcal{C}(\mathcal{H}_2).$ Suppose $\mathcal{H}_1$ is
derived equivalent to the module category of a hereditary algebra
$H$. Without loss the generality (compare Corollary 2.4.), we may assume that $\mathcal{H}_1=\mbox{mod}H.$
Then $\alpha (H)$ is a tilting object in
$\mathcal{C}(\mathcal{H}_2).$ By Proposition 2.6., $\alpha (H)$ is
induced by a tilting object $T'$ in $\mathcal{H}_2'$, derived
equivalent to $\mathcal{H}_2$. Then we have that
$$\begin{array}{lll}H&\cong
&\mbox{End}_{\mathcal{C}(\mathcal{H}_1)}(H) \ \ \ (\mbox{by
 Proposition 3.1}) \\
&\cong & \mbox{End}_{\mathcal{C}(\mathcal{H}_2)}(\alpha (H))\\
&\cong & \mbox{End}_{\mathcal{C}(\mathcal{H}'_2)}(T')\  \ \
(\mbox{by Corollary 2.4}) \\
 &\cong & \mbox{End}_{\mathcal{H}'_2}(T')\ltimes
D\mbox{Hom}_{\mathcal{H}'_2}(T', \tau ^2 T') \ \ (\mbox{by
 Proposition 3.1}).\end{array}$$
This shows $\mbox{End}_{\mathcal{H}'_2}(T')\ltimes
D\mbox{Hom}_{\mathcal{H}'_2}(T', \tau ^2 T')$ is hereditary, hence
$\mbox{End}_{\mathcal{H}'_2}(T')$ is hereditary [FGR], which is
denoted by $A$. It is not difficult to see that
$D\mbox{Hom}_{\mathcal{H}'_2}(T', \tau ^2 T')\cong
\mbox{Ext}^2_{A}(DA,A),$ hence it is zero.   Then
$\mbox{End}_{\mathcal{H}'_2}(T')\cong H.$ Therefore
$D^b(\mathcal{H}_i)$ is triangle equivalent to $D^b(H)$ for $i=1,
2$. The proof is finished.

\medskip

 If we restrict our attention to equivalences induced by tilting
 objects in $\mathcal{H},$ we can get some applications of Proposition 2.3.
 and Theorem 4.1. For any tilting object $T$ in $\mathcal{H}$,
Hom$_{\mathcal{H}}(T,-)$ induces a triangle equivalence from
$\mathcal{C(H)}$ to the cluster category $\mathcal{C}(A)$, where
$A$ is the quasi-tilted algebra of $T$. When $\mathcal{H}$ is the
category of finitely generated left modules over a hereditary
algebra, if we set
 $T$ to be a Bernstein-Gel'fand-Ponomarev titling module [BGP] or an
  APR-tilting
module [APR], then we get a triangle equivalence $R(T)_{\tau
^{-1}[1]}$ between the
 corresponding cluster categories. This triangle equivalence
 provides a realization of "truncate
reflection functors" in [FZ4][MRZ], if we identify the cluster
categories in each side with the set of almost positive roots of
corresponding Kac-Moody Lie algebras. We will explain simply in the
following, for details, we refer to [Z1] [MRZ].
\medskip

Let $(\G,\b)$ be a valued graph without cycles, $\Omega$ an
orientation.
   For any  vertex $k\in \G, $ we can define a new orientation
  $s_k\Omega$ of $(\G,\b)$ by reversing the direction of arrows
  along all edges containing $k$. A vertex $k\in \G$ is said to be
  a sink (or a source) with respect to $\Omega$ if there are no arrows
  starting (or ending) at vertex $k$.

  \medskip

  Let $K$ be a field and $(\G,\b, \Omega)$ a valued quiver.
   Let
  $\m=(F_i, {}_iM_j)_{i,j\in \G}$ be a reduced $K-$species of
  type $ \Omega; $ that
  is, for all $i, j \in \G$, $_iM_j$ is an $F_i-F_j-$bimodule,
  where $F_i$ and $F_j$
   are finite extensions of $K$  and dim$(_{i}M_{j})_{F_j}=d_{ij}$
   and dim$_{K}F_i=
   \varepsilon_i$. An $K-$representation $V=(V_i,{}_j\varphi_i)$ of
   $\m$ consists of $F_i-$ vector space $V_{i}, i\in \G$, and of an
   $F_j-$linear map $_j\varphi_i: V_i\otimes {}_iM_j\rightarrow V_j$
   for each arrow $i\rightarrow j$. Such representation is called
   finite dimensional if $\sum _{i\in \Gamma }\mbox{dim}_{K}V_i<\infty.$
    The
   category of finite-dimensional representations of $\m$ over
   $K$ is denoted by rep$(\m , \G, \Omega)$.

   \medskip

Now we fix an $K-$species $\m$ of a given valued quiver $(\G,\b,
\Omega)$. Given a sink, or a source $k$ of the quiver
    $(\G,\b,
\Omega)$,  we are going to recall the reflection functor
    $S^{\pm}_k$:

    $$S^+_k :\  \mbox{rep}(\m ,\G,  \Omega) \longrightarrow  \mbox{rep}(\m
    , \G, s_k\Omega),\  \mbox{    if } \ k  \ \mbox{ is a sink,}$$ or
$$S^-_k :\  \mbox{rep}(\m ,\G,  \Omega) \longrightarrow  \mbox{rep}(\m
,\G, s_k\Omega), \ \mbox{    if } \ k  \ \mbox{ is a source }.$$

\medskip

 We assume $k$ is a sink. For any representation $V=(V_i, \phi _{\alpha})$
  of
$(\m , \G, \Omega)$, the image of it under $S^+_k$ is by definition,
$S^+_kV=(W_i,{}_j\psi_i),$ a representation of $(\m , \G,
s_k\Omega)$, where $W_i=V_i,$ if $i\not=k;$ and $W_k$ is the kernel
in the diagram:
$$\begin{array}{lcccccccc}
(*)&&0&\longrightarrow& W_k&\s{({}_j\chi
_{k})_{j}}{\longrightarrow}& \oplus {}_{j\in \G} V_j\otimes {}_jM_k
&\s{({}_k\phi _{j})_j}{\longrightarrow}&V_k
\end{array}$$
${}_j\psi _i={}_j\phi _{i}$ and ${}_j\psi
 _{k}={}_j\bar{\chi}_{k}: W_k\otimes {}_kM_j\rightarrow X_j,$
 where ${}_j\bar{\chi}_{k}$ corresponds to ${}_j\chi _{k}$ under
 the isomorphism Hom$_{F_j}(W_{k}\otimes
 {}_kM_j,V_j)\approx\mbox{Hom}_{F_i}(W_k, V_j\otimes
 {}_jM_i ).$

If \textbf{$\alpha$}$=(\alpha _{i}): V\rightarrow V'$ is a morphism
in rep$(\m ,\G,  \Omega)$, then
$S^+_k$\textbf{$\alpha$}$=\beta=(\beta _i)$, where $\beta _i =\alpha
_i $ for $i\not=k$ and $\beta _k: W_k\rightarrow W_k'$ as the
restriction of $\oplus _{j\in \G}(\alpha _j\otimes 1)$ given in the
following commutative diagram:

\[ \begin{CD}
0@>>>W_k @>(_j\chi_{k})_{j}>>\oplus_{j\in \G} V_j\otimes
{}_jM_k @>(_k\phi_j)_j>>V_k\\
@VVV@VV\beta _k V @VV\oplus _j(\alpha _j\otimes 1)V @VV\alpha _k V \\
0@>>>W_k'@>(_j\chi_{k}')_{j}>>\oplus_{j\in \G} V_j'\otimes
{}_jM_k@>(_k\phi_j')_j>>V_k'
    \end{CD} \]

\medskip

If $k$ is source, the definition of $S^-_kV$ is dual to that of
$S^+_kV$, we omit it and refer to [DR].

\medskip

 In the rest of the section, we denote by $\cal{H}$ the category
 rep$(\m, \G, \Omega)$ and by
$\cal{H}'$ the category rep$(\m, \G, s_k\Omega)$, where $k$ is a
sink (or source) of $(\G,\b, \Omega)$. The root categories
$D^b(\mathcal{H})/[2]$, $D^b(\mathcal{H}')/[2]$ are denoted by
$\mathcal{R}(\Omega)$ and $\mathcal{R}(s_k\Omega)$ respectively. The
cluster categories $D^b(\mathcal{H})/F$, $D^b(\mathcal{H}')/F$ are
denoted by $\mathcal{C}(\Omega)$ and $\mathcal{C}(s_k\Omega)$
respectively.

\medskip

Let $P_i$ (or $P_i'$) be the projective indecomposable
representation $\mathcal{H}$ (resp. $\mathcal{H}'$) corresponding to
the vertex $i\in \G_0$, and  $E_k$ (or $E_k'$) the simple
 representation of $\mathcal{H}$ (resp. $\mathcal{H}'$)
 corresponding to the vertex $k$. We denote by $H$ (or $H'$) the
  tensor algebra
of $(\m, \G, \Omega)$ ($(\m, \G, s_k\Omega)$, resp.). Note that if
$k$ is a sink, then $P_k=E_k$.
 \medskip

  Let $T=\oplus _{i\in \G-{k}}P_i \oplus \tau ^{-1}P_k$. Suppose
$k$ is a sink, then $T$ is a tilting module in
$\mathcal{H}=$rep$(\m,\G, \Omega)$ which is called BGP-tilting
module (or APR-tilting). $S^+_k=\mbox{Hom}(T,-)$ as functors.

\medskip

The following lemma is proved in [Z1], for completeness, we give a
proof here.

\medskip

 \textbf{Lemma 4.2. }{\it Let $k$ be a sink (or a source)
  of a valued quiver $(\G,\b, \Omega).$ Then the BGP-reflection
   functor induces
  a triangle equivalence $R_{\tau ^{-1}[1]}(S^+_k)$
  (resp.,$R_{\tau ^{-1}[1]}(S^-_k)$) from
  $\mathcal{C}(\Omega)$ to
  $\mathcal{C}(s_k\Omega).$ Moreover $R_{\tau ^{-1}[1]}(S^+_k)
  (\widetilde{E_k})=
  \widetilde{P_k'}[1],$
   $R_{\tau ^{-1}[1]}(S^+_k)(\widetilde{P_k}[1])
   =\widetilde{E_k'},$ and for $j\neq k$,
    $R_{\tau ^{-1}[1]}(S^+_k)(\widetilde{P_j}[1])=\widetilde{P_j'}[1],$
     for indecomposable
    non-projective $H-$module $X,$  $R_{\tau ^{-1}[1]}(S^+_k)
    (\widetilde{X})=
    \widetilde{S^+_k(X)}.$ }
\medskip

\textbf{Proof.} From Corollary 2.5., we have the triangle
equivalent functor $R_{\tau ^{-1}[1]}(S^+_k)$ from the cluster
category $\mathcal{C}(\Omega)$ to $\mathcal{C}(s_k\Omega).$  Now
we prove that
\newline $R_{\tau ^{-1}[1]}(S^+_k)(\widetilde{E_k})=
  \widetilde{P_k'}[1].$ From [APR], we have AR-sequence
  $ (*):\ 0\rightarrow
E_k\rightarrow
  X\rightarrow\tau^{-1}E_k\rightarrow 0$ in $H-$mod with
   $X$ and $\tau^{-1}E_k$ being modules without $E_k$ as direct
   summands. Since $S^+_k$ is
    left exact functor, we have the exact sequence
 $0\rightarrow
  S^+_k(X)\rightarrow S^+_k(\tau^{-1}E_k)$ in $H'-$mod, in which
  the cokernal of the injective map is $E_k'.$
Regarded as the stalk complex
   of degree $0$, $E_k^{\bullet}$ is isomorphic to the complex:
  $\cdots\rightarrow 0\rightarrow X\rightarrow \tau^{-1}E_k\rightarrow 0
  \rightarrow\cdots$
  in $D^b(\mathcal{H}).$ By applying $S^+_k$ to the complex above,
   we have
  that $S^+_k(E_k^{\bullet})=\cdots\rightarrow 0\rightarrow S^+_k(X)
  \rightarrow S^+_k( \tau^{-1}E_k)\rightarrow0\rightarrow\cdots.$
   It follows that the complex $\cdots\rightarrow 0\rightarrow
S^+_k(X)\rightarrow S^+_k(\tau ^{-1}E_k)\rightarrow 0
\rightarrow\cdots$ is quasi-isomorphic to the stalk complex
$E_k'^{\bullet}[-1]$  of degree $-1$. It follows
 $R(S^+_k)(\widetilde{E_k})=\widetilde{E_k'}[-1].$ Since
$\tau \widetilde{P_k'} =\widetilde{E_k'}[-1], $ $R_{\tau
^{-1}[1]}(S^+_k)(\widetilde{E_k})=\tau \widetilde{P_k'} =
\widetilde{F^{-1}(P_k')}[1]=\widetilde{P_k'}[1].$
  The proof for others is easy: in derived category $D^b(\mathcal{H})$,
  we have that $S^+_k(P_i)=P_i'$ for any $i\neq k$,  $S^+_k(E_k[1])=E_k'.$
      It
follows that $R_{\tau
^{-1}[1]}(S^+_k)(\widetilde{P_i})=\widetilde{P_i'}$ for any $i\neq
k$ and $R_{\tau
^{-1}[1]}(S^+_k)(\widetilde{E_k}[1])=\widetilde{E_k'}$. Let $X\in
\mbox{ind}H$ be non-projective representation, $R_{\tau
^{-1}[1]}(S^+_k)(\widetilde{X})= \widetilde{R(S^+_k)(X)}=
    \widetilde{S^+_k(X)}.$ The proof is finished.

\medskip

 Let $(\G,\b, \Omega)$ be a valued quiver.
 We denote by $\Phi$ the set of roots of the
corresponding Kac-Moody Lie algebra. Let $\Phi_{\geq -1}$ denote the
set of almost positive roots, i.e.\ the positive roots together with
the negatives of the simple roots. When $\G$ is of Dynkin type, the
cluster variables of type $\G$ are in 1--1 correspondence with the
elements of $\Phi_{\geq -1}$ (compare [FZ3],[FZ4]).
 We define the map $\gamma_{\Omega}$ from
$\ind\mathcal{C}(\Omega)$ and $\Phi_{\geq -1}$ as follows (compare
[BMRRT]) : Let $X\in \mbox{ind}(\mbox{mod}H\vee H[1]).$
$$\gamma_{\Omega}(\widetilde{X})=\{
\begin{array}{lrl}\mathbf{dim}X & \mbox{ if } & X\in \mbox{ind}H;\\
&&\\
-\mathbf{dim}E_i& \mbox{ if } &X=P_i[1],\end{array}$$ where
$\mathbf{dim}X$ denotes the dimension vector of representation $X$.
When $\G$ is of Dynkin type, $\Phi_{\geq -1}$ is a bijection which
sends basic tilting objects to clusters in $\Phi_{\geq -1}$ (compare
[Z1]).

\medskip

  Let $s_i$ be the Coxeter generator of
Weyl group of $\Phi$ corresponding to $i\in \G $. We recall from
[FZ3][FZ4] that the "truncated reflection" $ \sigma_i$ of
$\Phi_{\geq -1}$ is defined as follows:
$$  \sigma_i(\alpha)=\left\{ \begin{array}{ll} \alpha &
\alpha=-\alpha_j,\ j\not=i \\ s_i(\alpha) & \mbox{otherwise.}
 \end{array}\right.$$

\medskip

 By using Lemma 4.2., one gets the following commutative diagram which
 explains the $R_{\tau ^{-1}[1]}(S_k^{+})$ is the realization of
 "truncate reflections" in [FZ3-4].

 \medskip

 \textbf{Proposition 4.3. }{\it Let $k$
be a sink (or a source) of a valued quiver $(\G,\b, \Omega)$. Then
we have the commutative diagram:
\[ \begin{CD}\mbox{ind}\mathcal{C}(\Omega) @>R_{\tau ^{-1}[1]}(S^+_k)
>(R_{\tau ^{-1}[1]}(S^-_k),resp.)> \mbox{ind}\mathcal{C}(s_k\Omega)
\\
@V\gamma _{\Omega} VV  @VV\gamma  _{s_k\Omega}V  \\
\Phi_{\geq -1} @>\sigma _k>> \Phi_{\geq -1}
\end{CD} \]}

\medskip

\medskip

\textbf{Remark 4.4.} If $(\G,\b, \Omega)$ is a simply-laced quiver
of Dynkin type, say $Q$, for a sink or source
 $k$, there are functors
 $\Sigma ^+_k$ and $\Sigma ^-_k$ respectively defined in [MRZ]
 which give a realization
 of $\sigma _k$ via "decorated" quiver representation. We remark
 that the
functors $\Sigma ^+_k$ and $\Sigma ^-_k$
 defined in [MRZ] are not equivalent.
 In this case, the functors $R_{\tau ^{-1}[1]}(S^+_k)$ also satisfy
  some community with functors $\Sigma ^+_k$
  in the following diagram: in the
 diagram, rep$\widetilde{Q}$ denotes the category of decorated
 representations of $Q$ and $\mathbf{sdim}(M)$ its signed dimension
 vector (refer
 [MRZ]). For $\Psi _Q$, we refer Section 4 in [BMRRT].

\[ \begin{CD}
\mbox{ind}\mathcal{C}(Q) @>R_{\tau ^{-1}[1]}(S^+_k)
>(R_{\tau ^{-1}[1]}(S^-_k),resp.)> \mbox{ind}\mathcal{C}(s_kQ)\\
@V\Psi_{Q} VV  @VV\Psi _{s_kQ}V\\
\mbox{indrep}\widetilde{Q }@>\Sigma ^+ _k>(\Sigma ^- _k,resp.)>
 \mbox{indrep}\widetilde{s_{k}Q} \\
@V\mathbf{sdim} VV  @VV\mathbf{sdim}V\\
\Phi_{\geq -1} @>\sigma _k>> \Phi_{\geq -1}
\end{CD} \]
\medskip

 We now apply Corollary 2.5. to the root categories $\mathcal{R}
 (\Omega)$, where $(\m, \G, \Omega)$
 is a valued quiver with species $\m$. The Grothendieck group $K_0(
 \mathcal{R}(\Omega))$ is $\mathbf{Z}^n$,
  where $n$ is the number of vertices of $\G$. For $M\in K_0(\mathcal{R}
  (\Omega))$, we denote by $\mathbf{dim}M$
  the canonical image of $M$ in $K_0(\mathcal{R}(\Omega)).$  It is easy
  to see that $\mathbf{dim}M[1]=-\mathbf{dim}M.$
  It follows from Kac's theorem [Ka] [DX] that we have a map
  $\mathbf{dim}:\
  \mbox{ind}\mathcal{R}(\Omega)\longrightarrow \Phi (\G)$, which
  is surjective. In case $\G$ is of Dynkin type,
  the map $\mathbf{dim}$
  is bijective. From Corollaries 2.4., 2.5., we get a similar
  commutative diagram to Proposition 4.3 (compare [XZZ]).

\medskip

\textbf{Proposition 4.5. }{\it Let $k$ be a sink (or a source)
  of a valued quiver $(\G,\b, \Omega)$ (of any type). Then the
  BGP-reflection functor induces
   a triangle equivalence $R_{[2]}(S^+_k)$
  ($R_{[2]}(S^-_k)$, resp.) from
  $\mathcal{R}(\Omega)$ to
  $\mathcal{R}(s_k\Omega),$ with $R_{[2]}(S^+_k)(\widetilde{E_k})=
  \widetilde{E_k'}[1].$ Moreover, we have the commutative diagram:
\[ \begin{CD}
\mbox{ind}\mathcal{R}(\Omega) @>R_{[2]}(S^+_k)
>(R_{[2]}(S^-_k),resp.)> \mbox{ind}\mathcal{R}(s_k\Omega)
\\
@V\mathbf{dim} VV  @VV\mathbf{dim}V  \\
\Phi @>s _k>> \Phi
\end{CD} \]}

\begin{center}
\textbf {ACKNOWLEDGMENTS.}\end{center} This work was completed and
revised  when the author was visiting University of Leicester
supported by Asia Link Programme "Algebras and Representations in
China and Europe" and was visiting Universitaet Paderborn and
Universitaet zu Koeln respectively.
 He would like to thank Professors Steffen Koenig and Henning Krause for valuable
 suggestions and
  comments and to them and the members of Algebras group
 there for their warm hospitality during his visits. He thanks Professors B. Keller, H. Krause,
 I.Reiten for helpful suggestions.   The author is grateful to the
 referee for a number of helpful comments and valuable suggestions.
%\newpage

\begin{center}

\end{center}

\medskip
\end{document}